\documentclass[12pt]{article}

\usepackage{amsmath, epsfig, cite}
\usepackage{amssymb}
\usepackage{amsfonts}
\usepackage{latexsym}
\usepackage{float}
\usepackage{color}
\usepackage{booktabs}
\usepackage{graphicx}
\usepackage{caption}
\usepackage[colorlinks=true,linkcolor=blue,citecolor=blue,urlcolor=blue]{hyperref}

\newtheorem{thm}{Theorem}[section]

\newtheorem{lem}[thm]{Lemma}

\newtheorem{ex}[thm]{Example}
\newcommand{\pf}{\noindent{\it Proof.} }

\numberwithin{equation}{section}

\newcommand{\qed}{{\hfill$\square$}\medskip}
\newcommand{\qbinom}[3]{\left[\begin{matrix}#1\\#2\end{matrix}\right]_{#3}}

\begin{document}

\begin{center}
{\Large\bf Truncated MacMahon-type $q$-series in arithmetic progressions and Chebyshev expansions}
\end{center}

\vskip 2mm \centerline{Ji-Cai Liu}
\begin{center}
{\footnotesize Department of Mathematics, Wenzhou University, Wenzhou 325035, PR China\\[5pt]
{\tt jcliu2016@gmail.com} \\[10pt]
}
\end{center}

\vskip 0.7cm \noindent{\bf Abstract.}
We study three finite families of MacMahon-type $q$-series.  The first is the arithmetic-progression family
$B_{k,m,N,r}^{\pm}(q)$, which contains the usual truncated MacMahon series and their odd-part analogues as special cases.  The other two families, $A_{k,m,a}(q)$ and $C_{k,m,a}(q)$, replace the squared linear factors by quadratic factors and are naturally controlled by Chebyshev polynomials.  We prove a double $q$-binomial expansion for $B_{k,m,N,r}^{\pm}(q)$ by a weighted elementary-symmetric-function argument; this proves and extends Merca's Conjecture 7 on truncated MacMahon series. We also establish finite Chebyshev polynomial expansions for $A_{k,m,a}(q)$ and $C_{k,m,a}(q)$, and we give weighted partition proofs of the $a=-2$ identities for both the ordinary family $A_{k,m,a}(q)$ and the odd family $C_{k,m,a}(q)$.

\vskip 3mm \noindent\parbox{\textwidth}{\raggedright {\it Keywords}: MacMahon's $q$-series; partitions; overpartitions; $q$-binomial coefficients; Chebyshev polynomials; finite Jacobi triple product}

\vskip 2mm
\noindent{\it MR Subject Classifications}: 05A17, 05A19, 11P81, 11F11

\section{Introduction}

A partition of a non-negative integer $n$ is a non-increasing sequence of positive integers whose sum is $n$.  Equivalently, one may write
\[
        n=t_1+2t_2+\cdots+nt_n,
\]
where $t_i$ is the multiplicity of the part $i$.  For standard background on partitions, Gaussian polynomials and $q$-series notation, see Andrews' text \cite{andrews-b-1998}.  This elementary notation already suggests a bridge between partitions and divisor functions: partitions of $n$ into identical parts are counted by the divisors of $n$, and weighting a part size by its multiplicity leads to the classical sum-of-divisors function
\[
        \sigma_1(n)=\sum_{d\mid n}d,
        \qquad
        \sum_{n\ge 1}\sigma_1(n)q^n=\sum_{s\ge 1}\frac{q^s}{(1-q^s)^2}.
\]

MacMahon's 1920 work on divisors and partitions \cite{macmahon-plms-1920} systematically developed this point of view.  For a fixed number of distinct part sizes, he considered generating functions whose coefficients are sums of products of multiplicities.  In modern notation, one of his main families is
\[
   A_k^+(q)=\sum_{1\le s_1<\cdots<s_k}
        \frac{q^{s_1+\cdots+s_k}}{(1-q^{s_1})^2\cdots(1-q^{s_k})^2},
\]
and the corresponding odd-part family is
\[
   C_k^+(q)=\sum_{1\le s_1<\cdots<s_k}
        \frac{q^{2s_1+\cdots+2s_k-k}}
        {(1-q^{2s_1-1})^2\cdots(1-q^{2s_k-1})^2}.
\]
The coefficient of $q^n$ in $A_k^+(q)$ is the sum of products of the multiplicities over all partitions of $n$ with exactly $k$ distinct part sizes; $C_k^+(q)$ has the analogous interpretation for partitions using exactly $k$ distinct odd part sizes.

Andrews and Rose \cite{ar-jram-2013} connected these series with Chebyshev polynomials and the Jacobi triple product.  A particularly useful consequence of their work is the pair of single-sum expansions
\begin{align}
A_k^+(q)
&=\frac{1}{(q;q)_\infty^3}
  \sum_{n=k}^{\infty}(-1)^{n-k}\frac{2n+1}{2k+1}\binom{n+k}{2k}q^{n(n+1)/2},
\label{eq:AR-A-intro}\\[5pt]
C_k^+(q)
&=\frac{(-q;q)_\infty}{(q;q)_\infty}
  \sum_{n=k}^{\infty}(-1)^{n-k}\frac{2n}{n+k}\binom{n+k}{2k}q^{n^2}.
\label{eq:AR-C-intro}
\end{align}
They also proved recurrence relations for these functions and showed that $A_k^+(q)$ is a linear combination of quasimodular forms of weights at most $2k$.  Rose \cite{rose-rnt-2015} further developed this quasimodularity viewpoint, and Bachmann \cite{bachmann-rnt-2024} related the odd family to multiple Eisenstein series and quasimodular forms on $\Gamma_0(2)$.

Several later papers refined the arithmetic and combinatorial meaning of MacMahon's functions.  Amdeberhan, Andrews and Tauraso \cite{aat-rms-2024} introduced the quadratic-denominator analogues
\[
   U_t(a,q)=\sum_{1\le n_1<\cdots<n_t}
      \prod_{j=1}^{t}\frac{q^{n_j}}{1+a q^{n_j}+q^{2n_j}}
      \qquad(a=0,\pm1,\pm2),
\]
which specialize to the ordinary MacMahon family when $a=-2$.  Amdeberhan, Ono and Singh \cite{aos-am-2024} made the quasimodularity of the ordinary family explicit in terms of Eisenstein series and obtained congruences.  Ono and Singh \cite{os-jcta-2024} proved the infinite tail identities
\begin{align}
\frac{1}{(q;q)_\infty^3}
&=q^{-k(k+1)/2}\sum_{m=k}^{\infty}\binom{2m+1}{m+k+1}A_m^+(q),
\label{eq:OS-A-intro}\\[5pt]
\frac{(-q;q)_\infty}{(q;q)_\infty}
&=q^{-k^2}\sum_{m=k}^{\infty}\binom{2m}{m+k}C_m^+(q),
\label{eq:OS-C-intro}
\end{align}
which express the generating functions for $3$-colored partitions and overpartitions, in the sense of Corteel and Lovejoy \cite{cl-tams-2004}, through tails of MacMahon's families.  Arithmetic aspects of related MacMahon-type functions, including the odd quadratic-denominator case, have also been studied by Sellers and Tauraso \cite{st-rj-2025,st-jcta-2026}; related convolution identities appear in \cite{xfy-rj-2025}.  Further recent work on nearby partition variants, traces and quasimodularity includes \cite{kms-rj-2025,matsusaka-rnt-2025}.

Finite forms of MacMahon's series and finite extensions of the Ono--Singh identities have also been studied by Xia \cite{xia-pems-2026}.  Merca \cite{merca-jcta-2025} developed truncated identities related to \eqref{eq:OS-A-intro} and \eqref{eq:OS-C-intro}.  To state Merca's results, recall that for $n\ge 0$,
\[
   (a;q)_n=(1-a)(1-aq)\cdots(1-aq^{n-1}),\qquad (a;q)_0=1,
\]
and, when $|q|<1$,
\[
   (a;q)_\infty=\prod_{j\ge 0}(1-aq^j).
\]
The $q$-binomial coefficient is
\[
   \qbinom{n}{k}{q}=\begin{cases}
   \displaystyle\frac{(q;q)_n}{(q;q)_k(q;q)_{n-k}},&0\le k\le n,\\[8pt]
   0,&\hbox{otherwise.}
   \end{cases}
\]
When another base is needed, we write it explicitly, for example
$\qbinom{n}{k}{q^2}$ or $\qbinom{n}{k}{q^N}$.
For $1\le k\le m$, define the finite ordinary and odd MacMahon series by
\begin{align*}
A_{k,m}^{\pm}(q)
&:=\sum_{1\le \lambda_1<\cdots<\lambda_k\le m}
\frac{q^{\lambda_1+\cdots+\lambda_k}}
{(1\mp q^{\lambda_1})^2\cdots(1\mp q^{\lambda_k})^2},\\[4pt]
C_{k,m}^{\pm}(q)
&:=\sum_{1\le \lambda_1<\cdots<\lambda_k\le m}
\frac{q^{2\lambda_1+\cdots+2\lambda_k-k}}
{(1\mp q^{2\lambda_1-1})^2\cdots(1\mp q^{2\lambda_k-1})^2}.
\end{align*}
Thus the superscript $+$ selects the minus sign in every denominator factor, whereas the superscript $-$ selects the plus sign.  With this notation, Merca proved
\begin{align*}
\sum_{j=k}^{m}\binom{2j+1}{j+k+1}A_{j,m}^+(q)
&=\frac{q^{k(k+1)/2}}{(q;q)_m^2}\qbinom{2m+1}{m+k+1}{q},\\[5pt]
\sum_{j=k}^{m}(\pm1)^{j-k}\binom{2j}{j+k}C_{j,m}^{\pm}(q)
&=\frac{q^{k^2}}{(\pm q;q^2)_m^2}\qbinom{2m}{m+k}{q^2}.
\end{align*}
He also proved finite analogues of \eqref{eq:AR-A-intro} and \eqref{eq:AR-C-intro}, namely
\begin{align}
(q;q)_m^2 A_{k,m}^+(q)
&=\sum_{n=k}^{m}(-1)^{n-k}\frac{2n+1}{2k+1}\binom{n+k}{2k}
  \qbinom{2m+1}{m+n+1}{q}q^{n(n+1)/2},
\label{eq:Merca-AR-A-intro}\\[5pt]
(q;q^2)_m^2 C_{k,m}^+(q)
&=\sum_{n=k}^{m}(-1)^{n-k}\frac{2n}{n+k}\binom{n+k}{2k}
  \qbinom{2m}{m+n}{q^2}q^{n^2}.
\label{eq:Merca-AR-C-intro}
\end{align}

Merca also formulated a double-sum conjecture for these same truncated series.  Conjecture 7 of \cite{merca-jcta-2025} asserts that, for positive integers $k\le m$,
\begin{align}
A_{k,m}^{\pm}(q)
&=\frac{1}{(\pm q;q)_m^2}
\sum_{i=0}^{m-k}\sum_{j=i+k}^{m}
(\mp1)^{j-i-k}\frac{2(j-i)}{j-i+k}\binom{j-i+k}{2k}
\qbinom{m}{i}{q}\qbinom{m}{j}{q}
q^{i(i+1)/2+j(j+1)/2},
\label{eq:Merca-conj-A-intro}\\[5pt]
C_{k,m}^{\pm}(q)
&=\frac{1}{(\pm q;q^2)_m^2}
\sum_{i=0}^{m-k}\sum_{j=i+k}^{m}
(\mp1)^{j-i-k}\frac{2(j-i)}{j-i+k}\binom{j-i+k}{2k}
\qbinom{m}{i}{q^2}\qbinom{m}{j}{q^2}
q^{i^2+j^2}.
\label{eq:Merca-conj-C-intro}
\end{align}
The first main theorem of the present paper, Theorem \ref{thm:B-main}, is a common generalization of these two conjectural identities: taking $(N,r)=(1,0)$ gives \eqref{eq:Merca-conj-A-intro}, while taking $(N,r)=(2,1)$ gives \eqref{eq:Merca-conj-C-intro}.  The elementary identity
\[
        \frac{2d}{d+k}\binom{d+k}{2k}
        =\frac{d}{k}\binom{d+k-1}{2k-1}
        \qquad(d\ge k)
\]
is the only change in the displayed coefficient.

The identities proved in this paper are finite extensions of these formulas.  The family $B_{k,m,N,r}^{\pm}(q)$ below extends the ordinary and odd squared-denominator cases to an arbitrary arithmetic progression $N\lambda-r$.  The families $A_{k,m,a}(q)$ and $C_{k,m,a}(q)$ extend \eqref{eq:Merca-AR-A-intro} and \eqref{eq:Merca-AR-C-intro} by replacing $(1-q^\ell)^2$ with the quadratic factor $1+a q^\ell+q^{2\ell}$.

Let $N,r$ be non-negative integers with $N>r$, and let $k,m$ be positive integers with $m\ge k$.  Define
\begin{align}
B_{k,m,N,r}^{\pm}(q)=
\sum_{1\le \lambda_1<\lambda_2<\cdots<\lambda_k\le m}
\frac{q^{N\lambda_1+\cdots+N\lambda_k-kr}}
{(1\mp q^{N\lambda_1-r})^2\cdots(1\mp q^{N\lambda_k-r})^2}.
\label{def:B}
\end{align}
The two finite families defined above are precisely the specializations
\[
        B_{k,m,1,0}^{\pm}(q)=A_{k,m}^{\pm}(q),\qquad
        B_{k,m,2,1}^{\pm}(q)=C_{k,m}^{\pm}(q).
\]

We also consider the quadratic-denominator analogues
\begin{align*}
A_{k,m,a}(q)=
\sum_{1\le \lambda_1<\cdots<\lambda_k\le m}
\frac{q^{\lambda_1+\cdots+\lambda_k}}
{\prod_{j=1}^{k}(1+a q^{\lambda_j}+q^{2\lambda_j})}.
\end{align*}
and
\begin{align*}
C_{k,m,a}(q)=
\sum_{1\le \lambda_1<\cdots<\lambda_k\le m}
\frac{q^{2\lambda_1+\cdots+2\lambda_k-k}}
{\prod_{j=1}^{k}(1+a q^{2\lambda_j-1}+q^{4\lambda_j-2})}.
\end{align*}
For the finite generating functions in Theorems \ref{thm:A-gen} and \ref{thm:C-gen}, we use the standard empty-product convention
\[
        A_{0,m,a}(q)=C_{0,m,a}(q)=1.
\]
The value $a=-2$ gives the squared denominators $(1-q^\ell)^2$ and therefore recovers the ``plus'' MacMahon-type cases.  The values $a=0,\pm1,\pm2$ are particularly natural because the roots of $1+ax+x^2$ are roots of unity or repeated roots.

The rest of the paper is arranged as follows.  Section 2 states the main results and explains their meaning.  Section 3 collects the auxiliary lemmas, including a direct proof of the binomial identity used in the double $q$-binomial expansion.  Section 4 proves the arithmetic-progression theorem for $B_{k,m,N,r}^{\pm}(q)$.  Section 5 proves the finite Chebyshev expansions for $A_{k,m,a}(q)$ and $C_{k,m,a}(q)$.  Sections 6 and 7 give weighted partition proofs of the $a=-2$ identities \eqref{eq:A-minus2} and \eqref{eq:C-minus2}, respectively.  The last section contains concluding remarks and further questions.

\section{Main results}

Our first result is a finite double $q$-binomial expansion for the arithmetic-progression family \eqref{def:B}.  It is a simultaneous extension of the corresponding finite formulas for the ordinary and odd MacMahon series.

\begin{thm}\label{thm:B-main}
Let $N,r$ be non-negative integers with $N>r$, and let $k,m$ be positive integers with $m\ge k$.  Then
\begin{align}
B_{k,m,N,r}^{\pm}(q)
&=\frac{1}{(\pm q^{N-r};q^N)_m^2}
\sum_{i=0}^{m-k}\sum_{j=i+k}^{m}(\mp 1)^{j-i-k}
\frac{j-i}{k}\binom{j-i+k-1}{2k-1}
\qbinom{m}{i}{q^N}\qbinom{m}{j}{q^N}
\notag\\[-2pt]
&\hspace{35mm}\times
q^{\frac{N i(i+1)}{2}-ir+\frac{N j(j+1)}{2}-jr}.
\label{eq:B-main}
\end{align}
\end{thm}

The theorem is finite and hence may be read formally as an identity of rational functions in $q$.  For $(N,r)=(1,0)$ and $(2,1)$ it proves Merca's Conjecture 7, because the coefficient in \eqref{eq:B-main} satisfies
\[
        \frac{j-i}{k}\binom{j-i+k-1}{2k-1}
        =\frac{2(j-i)}{j-i+k}\binom{j-i+k}{2k},
\]
which is precisely the coefficient appearing in \eqref{eq:Merca-conj-A-intro} and \eqref{eq:Merca-conj-C-intro}.  In the limit $m\to\infty$, it produces double-sum representations for the corresponding infinite MacMahon-type series in an arithmetic progression.

The next two theorems use Chebyshev polynomials.  Let $T_n$ be defined by
\[
        T_n(\cos\theta)=\cos(n\theta),
\]
and put
\[
        to_n(x)=\frac{T_{2n+1}(\sqrt{x})}{\sqrt{x}},\qquad
        te_n(x)=T_{2n}(\sqrt{x}).
\]
These are polynomials in $x$: $T_{2n+1}(u)/u$ and $T_{2n}(u)$ are even polynomials in $u$.  The square-root notation is therefore only a compact way to write these polynomials.

\begin{thm}\label{thm:A-gen}
For $m\ge 1$, we have
\begin{align}
\sum_{k=0}^{m}A_{k,m,a}(q)x^k
&=\prod_{i=1}^{m}\frac{1}{1+aq^i+q^{2i}}
\sum_{n=0}^{m}\qbinom{2m+1}{m+n+1}{q}
q^{n(n+1)/2}
 to_n\left(\frac{x+a+2}{4}\right).
\label{eq:A-gen}
\end{align}
Consequently, for $1\le k\le m$,
\begin{align}
A_{k,m,a}(q)
&=\prod_{i=1}^{m}\frac{1}{1+aq^i+q^{2i}}
\sum_{n=0}^{m}\qbinom{2m+1}{m+n+1}{q}
q^{n(n+1)/2}
[x^k]to_n\left(\frac{x+a+2}{4}\right).
\label{eq:A-coeff}
\end{align}
In particular,
\begin{align}
(q;q)_m^2 A_{k,m,-2}(q)
=\sum_{n=k}^{m}(-1)^{n-k}\frac{2n+1}{n+k+1}\binom{n+k+1}{2k+1}
\qbinom{2m+1}{m+n+1}{q}q^{n(n+1)/2}.
\label{eq:A-minus2}
\end{align}
\end{thm}

Theorem \ref{thm:A-gen} is a finite version of the Chebyshev expansion that underlies the Andrews-Rose formula for MacMahon's ordinary series.  The parameter $a$ moves the denominator from the repeated-root case $a=-2$ to other quadratic factors.

\begin{thm}\label{thm:C-gen}
For $m\ge 1$, we have
\begin{align}
\sum_{k=0}^{m}C_{k,m,a}(q)x^k
&=\prod_{i=1}^{m}\frac{1}{1+aq^{2i-1}+q^{4i-2}}\notag\\[-2pt]
&\times\left(
\qbinom{2m}{m}{q^2}
+2\sum_{n=1}^{m}\qbinom{2m}{m+n}{q^2}q^{n^2}
 te_n\left(\frac{x+a+2}{4}\right)
\right).
\label{eq:C-gen}
\end{align}
Consequently, for $1\le k\le m$,
\begin{align}
C_{k,m,a}(q)
&=2\prod_{i=1}^{m}\frac{1}{1+aq^{2i-1}+q^{4i-2}}
\sum_{n=1}^{m}\qbinom{2m}{m+n}{q^2}q^{n^2}
[x^k]te_n\left(\frac{x+a+2}{4}\right).
\label{eq:C-coeff}
\end{align}
In particular,
\begin{align}
(q;q^2)_m^2 C_{k,m,-2}(q)
=\sum_{n=k}^{m}\frac{2n(-1)^{n-k}}{n+k}\binom{n+k}{2k}
\qbinom{2m}{m+n}{q^2}q^{n^2}.
\label{eq:C-minus2}
\end{align}
\end{thm}

Formula \eqref{eq:C-minus2} is the odd analogue of \eqref{eq:A-minus2}.  Its importance is that it expresses the truncated odd MacMahon series through a single finite sum involving Gaussian polynomials.  Section 7 gives a weighted partition proof of this identity, in parallel with the finite $q$-binomial interpretation of the ordinary case in Section 6.

\section{Auxiliary lemmas}

We begin with the binomial identity that is used in the proof of Theorem \ref{thm:B-main}.  Although such identities can often be checked by computer-algebra methods such as the Wilf--Zeilberger framework \cite{pwz-b-1996}, we give a direct proof because the exact normalization is important in the symmetric-function argument below.

\begin{lem}\label{lem:binomial}
For non-negative integers $a$ and positive integers $k$, one has
\begin{align}
\sum_{j=\left\lceil (a+k)/2\right\rceil}^{a}
\frac{-a+2j}{k}\binom{-a+2j+k-1}{2k-1}\binom{a}{j}
=2^{a-k}\binom{a}{k}.
\label{eq:binomial-lemma}
\end{align}
\end{lem}

\pf
If $a<k$, then the sum is empty and both sides vanish.  Assume $a\ge k$ and put
\[
        P_k(X)=\frac{X}{k}\binom{X+k-1}{2k-1}.
\]
When a binomial coefficient has a possibly negative or formal upper argument, it is interpreted in the usual polynomial sense, equivalently by coefficient extraction from $(1+z)^Y$.  This convention is used only for finitely many coefficients.
The product form
\[
        P_k(X)=\frac{X^2}{k(2k-1)!}\prod_{h=1}^{k-1}(X^2-h^2)
\]
(with the empty product interpreted as $1$) shows that $P_k(X)$ is an even polynomial in $X$ and that $P_k(X)=0$ for integers $0\le X<k$.  We shall also use the elementary identity
\begin{align}
        P_k(X)=\binom{X+k}{2k}+\binom{X+k-1}{2k},
\label{eq:P-identity}
\end{align}
which follows after writing the two binomial coefficients over a common denominator.

The summand in \eqref{eq:binomial-lemma} is $P_k(2j-a)\binom{a}{j}$.  Because $P_k$ is even, the terms with indices $j$ and $a-j$ in the full sum
\[
        S_{a,k}=\sum_{j=0}^{a}P_k(2j-a)\binom{a}{j}
\]
are equal.  The possible middle term has $2j-a=0$ and hence is zero.  Moreover, any term with $0<2j-a<k$ is zero by the preceding paragraph.  Therefore the sum in \eqref{eq:binomial-lemma} is exactly one half of $S_{a,k}$.

Using \eqref{eq:P-identity} and coefficient extraction, we obtain
\begin{align*}
S_{a,k}
&=[z^{2k}]\sum_{j=0}^{a}\binom{a}{j}
\left((1+z)^{2j-a+k}+(1+z)^{2j-a+k-1}\right)\\
&=[z^{2k}](1+z)^{k-a-1}(2+z)
\left(1+(1+z)^2\right)^a.
\end{align*}
Since
\[
        1+(1+z)^2=2(1+z)\left(1+\frac{z^2}{2(1+z)}\right),
\]
we have
\begin{align*}
S_{a,k}
&=2^a\sum_{s=0}^{a}\binom{a}{s}2^{-s}
[z^{2k-2s}](2+z)(1+z)^{k-1-s}.
\end{align*}
For $s<k$, the polynomial $(2+z)(1+z)^{k-1-s}$ has degree $k-s$, which is smaller than the requested degree $2(k-s)$.  For $s>k$, the requested coefficient has negative degree and is zero.  Only $s=k$ remains, and its contribution is
\[
        2^a \binom{a}{k}2^{-k}[z^0](2+z)(1+z)^{-1}
        =2^{a-k+1}\binom{a}{k}.
\]
Thus $S_{a,k}=2^{a-k+1}\binom{a}{k}$, and the desired half-sum is \eqref{eq:binomial-lemma}.
\qed

The next lemma is the elementary symmetric function form of the same identity.

\begin{lem}\label{lem:symmetric}
Let $y_1,\ldots,y_m$ be indeterminates, let $e_s=e_s(y_1,\ldots,y_m)$ be the elementary symmetric function, and let $\eta\in\{1,-1\}$.  For $1\le k\le m$,
\begin{align}
\sum_{\substack{S\subseteq\{1,\ldots,m\}\\ |S|=k}}
\left(\prod_{s\in S}y_s\right)
\prod_{t\notin S}(1-\eta y_t)^2
=
\sum_{i=0}^{m-k}\sum_{j=i+k}^{m}(-\eta)^{j-i-k}
\frac{j-i}{k}\binom{j-i+k-1}{2k-1}e_i e_j.
\label{eq:symmetric}
\end{align}
\end{lem}

\pf
It is enough to compare the coefficient of each monomial
\[
        M=\prod_{u\in U}y_u^2\prod_{v\in V}y_v,
\]
where $U$ and $V$ are disjoint subsets of $\{1,\ldots,m\}$.  Put $t=|U|$, $\ell=|V|$, and $s=2t+\ell$.

If $\ell<k$, neither side can contain $M$, so both coefficients are zero.  We may therefore assume $\ell\ge k$.  On the left of \eqref{eq:symmetric}, the set $S$ must consist of exactly $k$ of the $\ell$ variables that occur to the first power in $M$.  The remaining $\ell-k$ variables in $V$ must then contribute the linear term $-2\eta y_v$ from $(1-\eta y_v)^2$, and each variable in $U$ contributes the quadratic term $y_u^2$.  Hence the coefficient of $M$ on the left is
\begin{align}
        (-2\eta)^{\ell-k}\binom{\ell}{k}.
\label{eq:left-coeff}
\end{align}

On the right, write the two summation indices temporarily as $\alpha$ and $\beta$.  The product $e_\alpha e_\beta$ can produce $M$ only when $\alpha+\beta=s$: every variable in $U$ must appear in both elementary-symmetric factors, while every variable in $V$ must appear in exactly one of them.  For fixed $\alpha,\beta$ with $\alpha+\beta=s$, the number of placements is
\[
        \binom{\ell}{\beta-t},
\]
because $\beta-t$ of the variables in $V$ are placed in the second factor.  Thus the coefficient of $M$ on the right is
\[
\sum_{\substack{\alpha,\beta\\ \alpha+\beta=s\\ \beta\ge \alpha+k}}
(-\eta)^{\beta-\alpha-k}
\frac{\beta-\alpha}{k}\binom{\beta-\alpha+k-1}{2k-1}
\binom{\ell}{\beta-t}.
\]
Set $h=\beta-t$.  Then $\beta-\alpha=2h-\ell$, and the condition $\beta\ge \alpha+k$ is equivalent to
$h\ge \left\lceil(\ell+k)/2\right\rceil$.  Also
$\beta-\alpha-k\equiv \ell-k\pmod 2$.  Therefore Lemma \ref{lem:binomial}, with $a=\ell$, gives
\[
\begin{aligned}
&\sum_{\substack{\alpha,\beta\\ \alpha+\beta=s\\ \beta\ge \alpha+k}}
(-\eta)^{\beta-\alpha-k}
\frac{\beta-\alpha}{k}\binom{\beta-\alpha+k-1}{2k-1}
\binom{\ell}{\beta-t}  \\
&\qquad =
(-\eta)^{\ell-k}
\sum_{h=\left\lceil(\ell+k)/2\right\rceil}^{\ell}
\frac{2h-\ell}{k}\binom{2h-\ell+k-1}{2k-1}\binom{\ell}{h}  \\
&\qquad =
(-\eta)^{\ell-k}2^{\ell-k}\binom{\ell}{k}.
\end{aligned}
\]
This is the same as \eqref{eq:left-coeff}.  The lemma follows.
\qed

We also need two finite Jacobi triple product forms.

\begin{lem}\label{lem:finite-A}
For $m\ge 1$,
\begin{align}
\sum_{n=0}^{m}\qbinom{2m+1}{m+n+1}{q}q^{n(n+1)/2}to_n\left(\frac{x}{4}\right)
=(q;q)_m^2\prod_{i=1}^{m}\left(1+\frac{xq^i}{(1-q^i)^2}\right).
\label{eq:finite-A}
\end{align}
\end{lem}

\pf
We first derive the finite triple product used here.  The finite $q$-binomial theorem gives
\[
        (-wq;q)_{2m+1}=\sum_{j=0}^{2m+1}\qbinom{2m+1}{j}{q}w^j q^{j(j+1)/2}.
\]
Put $j=m+n+1$ and $w=zq^{-m-1}$.  After multiplying by
$q^{m(m+1)/2}z^{-m-1}$ we obtain
\begin{align*}
\sum_{n=-m-1}^{m}\qbinom{2m+1}{m+n+1}{q}q^{n(n+1)/2}z^n
=(-zq;q)_m(-z^{-1};q)_{m+1}.
\end{align*}
Now set $z=e^{2i\theta}$ and multiply both sides by $e^{i\theta}$.  On the left, the terms with indices $n$ and $-n-1$ have the same Gaussian coefficient and the same power of $q$; pairing them gives
\[
        2\sum_{n=0}^{m}\qbinom{2m+1}{m+n+1}{q}q^{n(n+1)/2}T_{2n+1}(\cos\theta).
\]
On the right, the factor $(1+e^{-2i\theta})$ in $(-e^{-2i\theta};q)_{m+1}$ combines with $e^{i\theta}$ to give $2\cos\theta$, and the remaining factors give
\[
        2\cos\theta\prod_{i=1}^{m}(1+q^i e^{2i\theta})(1+q^i e^{-2i\theta}).
\]
Since
\[
(1+q^i e^{2i\theta})(1+q^i e^{-2i\theta})
=(1-q^i)^2+4q^i\cos^2\theta,
\]
the product equals
\[
        2\cos\theta\,(q;q)_m^2
        \prod_{i=1}^{m}\left(1+\frac{4q^i\cos^2\theta}{(1-q^i)^2}\right).
\]
Finally put $\cos\theta=\sqrt{x}/2$.  Since both sides are polynomials in $\cos^2\theta$ after division by $2\cos\theta$, this substitution is a formal polynomial substitution.  Dividing both sides by $2\cos\theta=\sqrt{x}$ changes the paired term
\[
        \frac{2T_{2n+1}(\sqrt{x}/2)}{\sqrt{x}}
\]
into
\[
        \frac{T_{2n+1}(\sqrt{x}/2)}{\sqrt{x}/2}=to_n(x/4),
\]
and proves \eqref{eq:finite-A}.
\qed

\begin{lem}\label{lem:finite-C}
For $m\ge 1$,
\begin{align}
2\sum_{n=1}^{m}\qbinom{2m}{m+n}{q^2}q^{n^2}te_n\left(\frac{x}{4}\right)
+\qbinom{2m}{m}{q^2}
=(q;q^2)_m^2\prod_{i=1}^{m}\left(1+\frac{xq^{2i-1}}{(1-q^{2i-1})^2}\right).
\label{eq:finite-C}
\end{align}
\end{lem}

\pf
We first record the finite triple product in the form needed below.  The finite $q$-binomial theorem with base $q^2$ gives
\[
(-zq^{1-2m};q^2)_{2m}
=\sum_{\nu=0}^{2m}\qbinom{2m}{\nu}{q^2}z^\nu q^{\nu^2-2m\nu}.
\]
Putting $\nu=m+n$ and multiplying by $z^{-m}q^{m^2}$ yields
\[
\sum_{n=-m}^{m}\qbinom{2m}{m+n}{q^2}q^{n^2}z^n
=z^{-m}q^{m^2}(-zq^{1-2m};q^2)_{2m}.
\]
In the product on the right, the exponents of $q$ are
\[
        1-2m,\,3-2m,\ldots,-1,1,\ldots,2m-1.
\]
Factoring $zq^e$ from the $m$ factors with negative exponent $e$ gives the cancelling factor
$z^m q^{-m^2}$ and leaves
\begin{align*}
\sum_{n=-m}^{m}\qbinom{2m}{m+n}{q^2}q^{n^2}z^n
=(-zq;q^2)_m(-z^{-1}q;q^2)_m.
\end{align*}
Set $z=e^{2i\theta}$.  Pairing the terms with indices $n$ and $-n$ gives
\[
\qbinom{2m}{m}{q^2}
+2\sum_{n=1}^{m}\qbinom{2m}{m+n}{q^2}q^{n^2}T_{2n}(\cos\theta)
=\prod_{i=1}^{m}(1+q^{2i-1}e^{2i\theta})(1+q^{2i-1}e^{-2i\theta}).
\]
The product on the right is
\[
(q;q^2)_m^2\prod_{i=1}^{m}
\left(1+\frac{4q^{2i-1}\cos^2\theta}{(1-q^{2i-1})^2}\right).
\]
Putting $\cos\theta=\sqrt{x}/2$ is again a formal polynomial substitution and gives $T_{2n}(\sqrt{x}/2)=te_n(x/4)$, which proves \eqref{eq:finite-C}.
\qed

Finally, the following coefficient evaluations are standard consequences of the explicit formula for Chebyshev polynomials.

\begin{lem}\label{lem:cheb-coeff}
For $0\le k\le n$,
\begin{align}
[x^k]to_n\left(\frac{x}{4}\right)
&=(-1)^{n-k}\frac{2n+1}{n+k+1}\binom{n+k+1}{2k+1},
\notag\\
[x^k]te_n\left(\frac{x}{4}\right)
&=(-1)^{n-k}\frac{n}{n+k}\binom{n+k}{2k}\qquad(n\ge 1).
\label{eq:te-coeff}
\end{align}
For $k>n$, both coefficients are zero.
\end{lem}

\pf
We use the standard formula
\[
T_N(x)=\frac{N}{2}\sum_{h=0}^{\lfloor N/2\rfloor}
(-1)^h\frac{(N-h-1)!}{h!(N-2h)!}(2x)^{N-2h}.
\]
For $to_n(x/4)$, take $N=2n+1$ and replace $x$ in the displayed formula by $\sqrt{x}/2$.  Since
\[
        to_n(x/4)=\frac{T_{2n+1}(\sqrt{x}/2)}{\sqrt{x}/2}
        =\frac{2T_{2n+1}(\sqrt{x}/2)}{\sqrt{x}},
\]
the term with $h=n-k$ contributes to $x^k$, and a direct simplification gives
\[
[x^k]to_n(x/4)
=(-1)^{n-k}\frac{2n+1}{n+k+1}\binom{n+k+1}{2k+1}.
\]
All other values of $h$ contribute different powers of $x$.  The proof of
\eqref{eq:te-coeff} is identical for $n\ge1$: put $N=2n$, replace $x$ by $\sqrt{x}/2$, and again set
$h=n-k$.  This gives
\[
[x^k]te_n(x/4)
=(-1)^{n-k}\frac{n}{n+k}\binom{n+k}{2k}.
\]
If $k>n$, no such value of $h$ occurs, and the coefficient is zero.
\qed

\section[Proof of Theorem \ref{thm:B-main}]{Proof of Theorem \ref{thm:B-main}}

Let
\[
        y_s=q^{Ns-r}\qquad(1\le s\le m).
\]
Then
\[
        e_i(y_1,\ldots,y_m)=q^{Ni(i+1)/2-ir}\qbinom{m}{i}{q^N}.
\]
Let $\eta=1$ when the superscript is $+$ and $\eta=-1$ when the superscript is $-$.  Then the denominator factor $1\mp y_s$ in \eqref{def:B} is $1-\eta y_s$, and
\[
        (\pm q^{N-r};q^N)_m=\prod_{s=1}^{m}(1-\eta y_s),
\]
where the upper signs are taken in the $+$ case and the lower signs in the $-$ case.  Multiplying \eqref{def:B} by this product squared gives
\begin{align*}
(\pm q^{N-r};q^N)_m^2 B_{k,m,N,r}^{\pm}(q)
&=\sum_{\substack{S\subseteq\{1,\ldots,m\}\\ |S|=k}}
\left(\prod_{s\in S}y_s\right)
\prod_{t\notin S}(1-\eta y_t)^2.
\end{align*}
The expression on the right is exactly the left side of Lemma \ref{lem:symmetric} with these choices of $y_s$ and $\eta$.  Applying that lemma yields
\begin{align*}
(\pm q^{N-r};q^N)_m^2 B_{k,m,N,r}^{\pm}(q)
&=\sum_{i=0}^{m-k}\sum_{j=i+k}^{m}(-\eta)^{j-i-k}
\frac{j-i}{k}\binom{j-i+k-1}{2k-1}e_i e_j.
\end{align*}
Since $-\eta$ is the sign denoted by $\mp 1$ in the statement of the theorem, substitution of the displayed formula for $e_i$ proves \eqref{eq:B-main}.
\qed

\section{Proofs of the Chebyshev expansions}

\noindent{\it Proof of Theorem \ref{thm:A-gen}.}
The left side is the elementary-symmetric generating function in the $m$ quantities
$q^i/(1+aq^i+q^{2i})$.  Hence
\begin{align*}
\sum_{k=0}^{m}A_{k,m,a}(q)x^k
&=\prod_{i=1}^{m}\left(1+\frac{xq^i}{1+aq^i+q^{2i}}\right)\\
&=\prod_{i=1}^{m}\frac{(1-q^i)^2}{1+aq^i+q^{2i}}
\prod_{i=1}^{m}\left(1+\frac{(x+a+2)q^i}{(1-q^i)^2}\right).
\end{align*}
Applying Lemma \ref{lem:finite-A} with $x$ replaced by $x+a+2$ gives \eqref{eq:A-gen}.  Extracting the coefficient of $x^k$ gives \eqref{eq:A-coeff}.  Finally, if $a=-2$, then the factor in front is $(q;q)_m^{-2}$ and Lemma \ref{lem:cheb-coeff} gives \eqref{eq:A-minus2}.
\qed

\noindent{\it Proof of Theorem \ref{thm:C-gen}.}
The same elementary-symmetric argument, now applied to the odd exponents $2i-1$, gives
\begin{align*}
\sum_{k=0}^{m}C_{k,m,a}(q)x^k
&=\prod_{i=1}^{m}\left(1+\frac{xq^{2i-1}}{1+aq^{2i-1}+q^{4i-2}}\right)\\
&=\prod_{i=1}^{m}\frac{(1-q^{2i-1})^2}{1+aq^{2i-1}+q^{4i-2}}
\prod_{i=1}^{m}\left(1+\frac{(x+a+2)q^{2i-1}}{(1-q^{2i-1})^2}\right).
\end{align*}
Applying Lemma \ref{lem:finite-C} with $x$ replaced by $x+a+2$ proves \eqref{eq:C-gen}.  Extracting the coefficient of $x^k$ with $k\ge 1$ gives \eqref{eq:C-coeff}; the central term \(\qbinom{2m}{m}{q^2}\) has no positive power of $x$.  With $a=-2$, Lemma \ref{lem:cheb-coeff} gives \eqref{eq:C-minus2}.
\qed

\section[A weighted partition proof of \eqref{eq:A-minus2}]{A weighted partition proof of \eqref{eq:A-minus2}}

We now give a coefficientwise proof of \eqref{eq:A-minus2}.  The proof is ``weighted'' because the same bounded partition may occur with a signed integer weight.  This is the finite analogue of the usual interpretation of MacMahon's series as sums of products of multiplicities.

Put $z_i=q^i$ for $1\le i\le m$.  From the definition of $A_{k,m,-2}(q)$ we have
\begin{align*}
(q;q)_m^2A_{k,m,-2}(q)
&=[x^k]\prod_{i=1}^{m}\left((1-z_i)^2+xz_i\right).
\end{align*}
Indeed, choosing $xz_i$ marks the part size $i$ as one of the $k$ distinguished singleton sizes; for an unmarked size one chooses one of the three terms $1,-2z_i,z_i^2$ coming from $(1-z_i)^2$.

Let $\mathcal P_{2,m,\ell}(N)$ be the set of partitions of $N$ in which every part is at most $m$, no part has multiplicity greater than $2$, and exactly $\ell$ part sizes occur with multiplicity $1$.  Write $P_{2,m,\ell}(N)=|\mathcal P_{2,m,\ell}(N)|$.  If a partition in $\mathcal P_{2,m,\ell}(N)$ is fixed, then $k$ of its $\ell$ singleton part sizes must be selected to contribute the factors $xz_i$.  The remaining $\ell-k$ singleton sizes contribute the factors $-2z_i$.  Therefore
\begin{align}
[q^N](q;q)_m^2A_{k,m,-2}(q)
=\sum_{\ell=0}^{m}(-2)^{\ell-k}\binom{\ell}{k}P_{2,m,\ell}(N),
\label{eq:A-left-weighted}
\end{align}
where, as usual, $\binom{\ell}{k}=0$ for $\ell<k$.

\begin{ex}
Take $m=6$ and $k=2$.  Consider the bounded partition
\[
        \pi=(6^2,5,4,3^2,2,1^2)
\]
of $31$.  The singleton part sizes are $2,4,5$, so $\ell=3$, while the part sizes $1,3,6$ occur twice.  In the product
\[
        \prod_{i=1}^{6}\left((1-q^i)^2+xq^i\right)
\]
this fixed partition contributes to $[x^2q^{31}]$ by choosing two of the three singleton sizes to supply the marked factors $xq^i$.  The remaining singleton size must then supply the term $-2q^i$ from $(1-q^i)^2$.  Thus the total contribution of this one partition is
\[
        \binom{3}{2}(-2)^{3-2}=-6,
\]
which is exactly the summand in \eqref{eq:A-left-weighted} for $\ell=3$.  This example illustrates why the proof is naturally weighted rather than purely enumerative.
\end{ex}

Next we interpret the Gaussian polynomial on the right side of \eqref{eq:A-minus2}.  Let $\mathcal Q_{m,n}(N)$ denote the set of partitions
$\lambda=(\lambda_1,\lambda_2,\ldots)$ of $N$ satisfying the following conditions:
\begin{itemize}
\item $\lambda_1\le m$ and the number of parts is at most $m+n+1$;
\item the first $n$ parts are strictly decreasing and the $n$-th part is larger than every remaining part, that is,
\[
      \lambda_1>\lambda_2>\cdots>\lambda_n>\lambda_{n+1},
\]
where missing parts are interpreted as $0$.
\end{itemize}
When $n=0$, the second condition is empty.  Let $Q_{m,n}(N)=|\mathcal Q_{m,n}(N)|$.  We claim that
\begin{align}
\sum_{N\ge0}Q_{m,n}(N)q^N
=q^{n(n+1)/2}\qbinom{2m+1}{m+n+1}{q}.
\label{eq:A-Gaussian-comb}
\end{align}
To see this, use the symmetry of Gaussian polynomials to write
\[       \qbinom{2m+1}{m+n+1}{q}=\qbinom{2m+1}{m-n}{q}.  \]
For $n=0$, \eqref{eq:A-Gaussian-comb} is just the standard interpretation of this Gaussian polynomial as partitions with at most $m+1$ parts, each at most $m$.  For $n\ge1$, the same Gaussian polynomial is the generating function for partitions $\mu$ with at most $m+n+1$ parts, each at most $m-n$.  Add the staircase $(n,n-1,\ldots,1)$ to the first $n$ parts of $\mu$ and leave all later parts unchanged.  This increases the weight by $n(n+1)/2$ and produces a partition in $\mathcal Q_{m,n}(N)$.

Conversely, subtract the same staircase from the first $n$ parts of a partition $\lambda\in\mathcal Q_{m,n}(N)$.  The strict inequalities among the first $n$ parts ensure that the resulting first $n$ parts are weakly decreasing, and the boundary inequality $\lambda_n>\lambda_{n+1}$ ensures that the new $n$-th part is at least the $(n+1)$-st.  The upper bound is also preserved: for $1\le i\le n$, strictness gives $\lambda_i\le \lambda_1-(i-1)$ and hence
\[
        \lambda_i-(n-i+1)\le \lambda_1-n\le m-n;
\]
for $i>n$, the boundary condition gives $\lambda_i\le \lambda_n-1\le \lambda_1-n\le m-n$.  Thus the inverse image is a partition $\mu$ with at most $m+n+1$ parts, each at most $m-n$.  Hence \eqref{eq:A-Gaussian-comb} follows.

\begin{ex}
Let $m=6$ and $n=3$.  The Gaussian polynomial
$\qbinom{13}{10}{q}=\qbinom{13}{3}{q}$ may be interpreted as partitions with at most $10$ parts, each at most $3$.  For instance,
\[
        \mu=(3,3,2,2,1,1,1)
\]
has weight $13$.  Adding the staircase $(3,2,1)$ to the first three parts gives
\[
        \lambda=(6,5,3,2,1,1,1),
\]
whose weight is $13+3+2+1=19$.  The first three parts satisfy
$6>5>3>2$, so $\lambda\in\mathcal Q_{6,3}(19)$.  Conversely, subtracting
$(3,2,1)$ from the first three parts of this $\lambda$ recovers $\mu$.  This example shows concretely how the factor $q^{n(n+1)/2}$ in \eqref{eq:A-Gaussian-comb} records the added staircase.
\end{ex}

Finally, Lemma \ref{lem:finite-A} says exactly that
\begin{align*}
\prod_{i=1}^{m}\left((1-q^i)^2+xq^i\right)
=\sum_{n=0}^{m}\qbinom{2m+1}{m+n+1}{q}q^{n(n+1)/2}to_n\left(\frac{x}{4}\right).
\end{align*}
Taking the coefficient of $x^k$ and using Lemma \ref{lem:cheb-coeff} gives \eqref{eq:A-minus2}.  Equivalently, for every $N\ge0$,
\begin{align}
\sum_{\ell=0}^{m}(-2)^{\ell-k}\binom{\ell}{k}P_{2,m,\ell}(N)
&=\sum_{n=k}^{m}(-1)^{n-k}\frac{2n+1}{n+k+1}\binom{n+k+1}{2k+1}Q_{m,n}(N).
\label{eq:A-coefficient-comb}
\end{align}
This coefficientwise identity is the desired combinatorial form of \eqref{eq:A-minus2}.

\section[A weighted partition proof of \eqref{eq:C-minus2}]{A weighted partition proof of \eqref{eq:C-minus2}}

The proof for the odd family is parallel to the preceding section, but the staircase and the allowed parts are odd.  Let $z_i=q^{2i-1}$ for $1\le i\le m$.  From the definition of $C_{k,m,-2}(q)$,
\begin{align*}
(q;q^2)_m^2C_{k,m,-2}(q)
=[x^k]\prod_{i=1}^{m}\left((1-z_i)^2+xz_i\right).
\end{align*}
Thus the left side enumerates partitions into odd parts not exceeding $2m-1$, again with multiplicities at most $2$, with $k$ distinguished singleton part sizes.

Let $\mathcal P^{\rm odd}_{2,m,\ell}(N)$ be the set of partitions of $N$ into odd parts not exceeding $2m-1$, with no part appearing more than twice, and with exactly $\ell$ part sizes appearing once.  Put $P^{\rm odd}_{2,m,\ell}(N)=|\mathcal P^{\rm odd}_{2,m,\ell}(N)|$.  The same marking argument as above gives
\begin{align*}
[q^N](q;q^2)_m^2C_{k,m,-2}(q)
=\sum_{\ell=0}^{m}(-2)^{\ell-k}\binom{\ell}{k}P^{\rm odd}_{2,m,\ell}(N).
\end{align*}

\begin{ex}
Take $m=7$ and $k=2$, so the allowed odd part sizes are
$1,3,5,7,9,11,13$.  The partition
\[
        \pi=(13^2,11,9,7^2,5,3^2,1)
\]
of $72$ has singleton odd part sizes $1,5,9,11$ and doubled odd part sizes $3,7,13$.  Hence $\ell=4$.  In
\[
        \prod_{i=1}^{7}\left((1-q^{2i-1})^2+xq^{2i-1}\right),
\]
the contribution of this fixed partition to $[x^2q^{72}]$ is
\[
        \binom{4}{2}(-2)^{4-2}=24.
\]
Thus the same signed singleton-counting mechanism used in the ordinary case survives unchanged, while only the permitted part sizes have changed.
\end{ex}

We now interpret the Gaussian polynomial in \eqref{eq:C-minus2}.  For $1\le n\le m$, let $\mathcal Q^{\rm odd}_{m,n}(N)$ be the set of partitions
$\lambda=(\lambda_1,\lambda_2,\ldots)$ of $N$ satisfying:
\begin{itemize}
\item all parts are at most $2m-1$, and the number of parts is at most $m+n$;
\item the first $n$ parts are distinct odd parts;
\item all remaining parts, if any, are even.
\end{itemize}
Let $Q^{\rm odd}_{m,n}(N)=|\mathcal Q^{\rm odd}_{m,n}(N)|$.  Then
\begin{align}
\sum_{N\ge0}Q^{\rm odd}_{m,n}(N)q^N
=q^{n^2}\qbinom{2m}{m+n}{q^2}.
\label{eq:C-Gaussian-comb}
\end{align}
Indeed, $\qbinom{2m}{m+n}{q^2}=\qbinom{2m}{m-n}{q^2}$ counts partitions $\mu$ into even parts, with at most $m+n$ parts and each part at most $2(m-n)$.  Add the odd staircase
\[
        (2n-1,2n-3,\ldots,3,1)
\]
to the first $n$ parts of $\mu$.  The added weight is $n^2$, the first $n$ resulting parts are distinct odd parts not exceeding $2m-1$, and all later parts remain even.

Conversely, subtract the same odd staircase from the first $n$ parts of a partition $\lambda\in\mathcal Q^{\rm odd}_{m,n}(N)$.  Since the first $n$ parts are distinct odd parts in a non-increasing partition, consecutive differences among them are at least $2$, so the first $n$ new parts are weakly decreasing even parts.  The $n$-th odd part is larger than the next even part, if such a part exists, so the boundary between the $n$-th and $(n+1)$-st parts is preserved.  The upper bound is preserved as follows: for $1\le i\le n$, the distinct odd condition gives $\lambda_i\le 2m-1-2(i-1)$, and hence
\[
        \lambda_i-(2n-2i+1)\le 2(m-n);
\]
for $i>n$, any later even part is at most $\lambda_n-1\le 2(m-n)$.  Therefore the inverse image is a partition into even parts with at most $m+n$ parts, each at most $2(m-n)$, and \eqref{eq:C-Gaussian-comb} follows.

The central term $\qbinom{2m}{m}{q^2}$ in Lemma \ref{lem:finite-C} corresponds to the case without an odd staircase; it has the separate standard interpretation as even partitions with at most $m$ parts, each at most $2m$.  It has no positive power of the marking variable $x$, so it does not enter the coefficient extraction for $k\ge1$ below.

\begin{ex}
Let $m=7$ and $n=3$.  Then
$\qbinom{14}{10}{q^2}=\qbinom{14}{4}{q^2}$ counts partitions into even parts with at most $10$ parts, each at most $8$.  Choose
\[
        \mu=(8,8,6,4,4,2,2),
\]
of weight $34$.  Adding the odd staircase $(5,3,1)$ to the first three parts gives
\[
        \lambda=(13,11,7,4,4,2,2),
\]
of weight $43=34+9$.  The first three parts are distinct odd parts, all later parts are even, and the largest part is $13=2m-1$.  Subtracting the same staircase from the first three parts recovers the original even partition $\mu$.  This illustrates the role of $q^{n^2}$ in \eqref{eq:C-Gaussian-comb}.
\end{ex}

Lemma \ref{lem:finite-C} gives
\begin{align}
\prod_{i=1}^{m}\left((1-q^{2i-1})^2+xq^{2i-1}\right)
=\qbinom{2m}{m}{q^2}
+2\sum_{n=1}^{m}\qbinom{2m}{m+n}{q^2}q^{n^2}te_n\left(\frac{x}{4}\right).
\label{eq:C-product-cheb}
\end{align}
For $k\ge1$, the central term has no contribution to $[x^k]$.  Taking $[x^k]$ in \eqref{eq:C-product-cheb} and using Lemma \ref{lem:cheb-coeff} proves \eqref{eq:C-minus2}.  Coefficientwise, this says that for every $N\ge0$,
\begin{align*}
\sum_{\ell=0}^{m}(-2)^{\ell-k}\binom{\ell}{k}P^{\rm odd}_{2,m,\ell}(N)
&=\sum_{n=k}^{m}\frac{2n(-1)^{n-k}}{n+k}\binom{n+k}{2k}Q^{\rm odd}_{m,n}(N).
\end{align*}
which is the desired weighted odd-part analogue of \eqref{eq:A-coefficient-comb}.

\section{Concluding remarks}

The three families considered here emphasize complementary aspects of MacMahon-type $q$-series.  The functions $B_{k,m,N,r}^{\pm}(q)$ show that the double $q$-binomial expansion is not tied to the full set of positive integers or to the odd integers alone; it is an elementary symmetric-function identity after the substitution $y_s=q^{Ns-r}$.  The functions $A_{k,m,a}(q)$ and $C_{k,m,a}(q)$ show that finite Chebyshev products naturally control quadratic denominators.  The case $a=-2$ links these products back to the classical squared-denominator MacMahon series.

Several questions remain natural.  First, the finite formulas in Theorems \ref{thm:A-gen} and \ref{thm:C-gen} should have further coefficientwise interpretations for $a=0,\pm1,\pm2$, especially when the quadratic denominator factors into roots of unity.  Second, the limiting forms as $m\to\infty$ connect the present identities to quasimodular forms and to the generating functions for $3$-colored partitions and overpartitions.  It would be interesting to make the passage from the finite weighted partition identities to infinite congruence families as explicit as possible.  Finally, the arithmetic-progression parameter pair $(N,r)$ in $B_{k,m,N,r}^{\pm}(q)$ suggests further refinements by residue classes, which may lead to modular or quasimodular behavior on congruence subgroups.

%

\end{document}